\documentclass[12pt]{article}
\usepackage{amsfonts, amsmath, amsthm, amssymb, fullpage}

\newtheorem{theorem}{Theorem}[section]
\numberwithin{equation}{theorem}
\newtheorem{lemma}[theorem]{Lemma}

\newtheorem{cor}[theorem]{Corollary}
\newtheorem{prop}[theorem]{Proposition}
\newtheorem{question}[theorem]{Question}

\theoremstyle{definition}
\newtheorem{example}[theorem]{Example}
\newtheorem{convention}[theorem]{Convention}
\newtheorem{defn}[theorem]{Definition}
\newtheorem{remark}[theorem]{Remark}

\def\CC{\mathbb{C}}
\def\FF{\mathbb{F}}
\def\QQ{\mathbb{Q}}

\def\ZZ{\mathbb{Z}}

\def\gothm{\mathfrak{m}}
\def\gotho{\mathfrak{o}}

\def\sep{\mathrm{sep}}

\DeclareMathOperator{\Aut}{Aut}
\DeclareMathOperator{\Di}{Di}
\DeclareMathOperator{\Gal}{Gal}
\DeclareMathOperator{\GL}{GL}
\DeclareMathOperator{\Ind}{Ind}
\DeclareMathOperator{\Norm}{Norm}
\DeclareMathOperator{\Stab}{Stab}

\newcounter{fixmectr}

\begin{document}

\title{Mass formulas for local Galois representations 
(with an appendix by Daniel Gulotta)}
\author{Kiran S. Kedlaya \\ Department of Mathematics \\ Massachusetts
Institute of Technology \\ 77 Massachusetts Avenue \\
Cambridge, MA 02139 \\
\texttt{kedlaya@math.mit.edu}}
\date{February 13, 2007}

\maketitle

\begin{abstract}
Bhargava has given a formula, derived from a formula of Serre,
computing a certain count of extensions of a local field,
weighted by conductor and by number of automorphisms.
We interpret this result as a counting formula for permutation
representations
of the absolute Galois group of the local field, then speculate on
variants of this formula in which the role of the symmetric group
is played by other groups. We prove an analogue of Bhargava's formula
for representations into a Weyl group in the $B_n$ series, which suggests
a possible link with integration on $p$-adic groups. We also obtain analogous
positive
results in odd residual characteristic, and negative results
in residual characteristic 2, for the $D_n$ series (in the appendix)
and the exceptional group $G_2$.
\end{abstract}

\section{Introduction}

Serre \cite[Th\'eor\`eme~2]{serre-mass} gave the following astonishing
``mass formula'' counting
totally ramified degree $n$ extensions of a local field $K$ with
residue field $\FF_q$: if $S_n$ is a set of representatives of
the isomorphism classes of such extensions of $K$, then
\begin{equation} \label{eq:serremass}
\sum_{L \in S_n} \frac{1}{w(L) q^{c(L)-n+1}} = 1,
\end{equation}
where $w(L)$ is the number of automorphisms of $L$ and
$c(L)$ is the discriminant exponent of $L$ over $K$. The automorphism
contribution is no surprise, as it invariably occurs in counting problems
of this ilk (essentially because of Burnside's formula); the distinguishing
feature of Serre's formula is the weighting by conductor.

In the context of deriving heuristics on the number of number fields
of given degree with discriminant bounded in a certain range
(consistent with the theorems of Davenport-Heilbronn 
\cite{davenport-heilbronn} in the cubic case
and Bhargava \cite{bhargava2, bhargava3} in the quartic and quintic 
cases; see Belabas's Seminaire Bourbaki notes \cite{belabas} for
an overview),
Bhargava \cite[Theorem~1]{bhargava} has derived from Serre's formula the 
following 
mass formula counting \'etale $K$-algebras of degree $n$. 
\begin{theorem}[Bhargava] \label{thm:abb}
Let $K$ be a local field with residue field $\FF_q$,
and let $\Sigma_n$ be a set of representatives for the isomorphism
classes of \'etale $K$-algebras of degree $n$. For $L \in \Sigma_n$,
let $w(L)$ be the number of automorphisms of $L$ and
let $c(L)$ be the discriminant exponent of $L$. Then one has
\begin{equation} \label{eq:bhargava}
\sum_{L \in \Sigma_n} \frac{1}{w(L) q^{c(L)}}
= \sum_{k=0}^{n} \frac{P(n,n-k)}{q^k},
\end{equation}
where $P(n,n-k)$ denotes the number of partitions of the integer $n$ into
exactly $n-k$ parts, or equivalently the number of partitions of $k$ into
at most $n-k$ parts. (Note: $P(n,0) = 1$ for $n=0$ and $0$ for $n>0$.)
\end{theorem}
The purpose of this paper is twofold. We first
reformulate Bhargava's formula as a
counting formula for permutation representations of the absolute Galois group
of a local field, and exhibit a straightforward deduction of the latter formula
from Serre's formula using standard techniques from combinatorics
(notably the Exponential Formula). We then pose some questions 
about possible mass formulas for other types of representations,
and establish affirmative and negative answers in some classes of
cases. More precisely, we ask (Question~\ref{question:weyl}) about 
representations into the Weyl group of
a semisimple Lie algebra, motivated by a potential link to integration
on $p$-adic Lie groups. Bhargava's formula answers
Question~\ref{question:weyl} affirmatively for the $A_n$ series; 
imitating the $A_n$ proof, we are able to resolve
Question~\ref{question:weyl} affirmatively for the $B_n$ series.
We also check the $G_2$ case by direct calculation; here we discover a
surprising negative answer to Question~\ref{question:weyl} for local
fields of residual characteristic 2, which renders any potential link to 
integration on $p$-adic groups even more mysterious.

In the appendix by Daniel Gulotta, Question~\ref{question:weyl} is answered
affirmatively for the $D_n$ series in the case of odd residual
characteristic. By machine calculation, it is also shown that
this affirmative answer cannot in general extend to residual characteristic 2.

\subsection*{Notational conventions}

In this paper, a \emph{local field} is a complete discretely valued
field (of either mixed or equal characteristics)
 with finite residue field. For $K$ a local field, let $\gotho_K$
denote the ring of integers of $K$, and let $\gothm_K$ denote the maximal ideal
of $\gotho_K$. 
For $L/K$ a finite separable extension of local fields, let $f(L/K)$ denote
the degree of the induced extension on residue fields,
and let $c(L/K)$ denote the discriminant exponent. For $K$ any field,
let $G_K = \Gal(K^{\sep}/K)$ denote the absolute Galois
group of $K$.

\subsection*{Acknowledgments}
Some of this material was presented at the workshop ``Rings of low rank''
in June 2006, hosted by the Lorentz Center in Leiden.
Thanks to Manjul Bhargava for providing a copy of his preprint
\cite{bhargava} and for additional helpful discussions,
and to Melanie Wood for comments and corrections on a previous draft.
Kedlaya was supported by NSF grant DMS-0400747, 
NSF CAREER grant DMS-0545904, and a Sloan Research Fellowship.
Gulotta was supported by MIT's Undergraduate Research Opportunities
Program.

\section{Total mass}

To begin with, we define a numerical invariant attached
to a local field and a finite linear group. First, let us fix 
notation for local fields.

\begin{defn}
Let $K$ be a local field, let $L$ be a finite Galois extension of $K$, 
and put $G = \Gal(L/K)$.
For $s \in G \setminus \{e\}$, let $i_G(s)$ be the largest integer
$n$ such that $v_L(x^s - x) \geq n$ for all $x \in \gotho_L$.
Define the Artin character $a_G: G \to \ZZ$ by
\[
a_G(s) = \begin{cases}
-f(L/K) i_G(s) & s \neq e \\
-\sum_{t \neq e} a_G(t) & s=e.
         \end{cases}
\]
By a theorem of Artin \cite[Theorem~VI.1]{serre},
$a_G$ is the character of a representation of $G$. Hence for any function
$\chi: G \to \CC$ which is the character of a complex representation of $G$,
the inner product
\[
\frac{1}{|G|} \sum_{s \in G} a_G(s) \chi(s)
\]
is a nonnegative integer, called the \emph{Artin conductor} of $\chi$.
For $\rho: G \to \GL_n(\CC)$ a representation with character $\chi$,
we write $c(\rho)$ for the conductor of $\chi$; note that
$c(\rho_1 \oplus \rho_2) = c(\rho_1) + c(\rho_2)$.
\end{defn}

\begin{defn} \label{D:total mass}
Let $n$ be a positive integer and let $\Gamma$ be a finite subgroup of 
$\GL_n(\CC)$. For $K$ a local field 
with residue field $\FF_q$, define the \emph{total mass} of the 
pair $(K,\Gamma)$, denoted $M(K, \Gamma)$, as follows.
Let $S_{K,\Gamma}$ be the set of continuous homomorphisms
$\rho: \Gal(K^{\sep}/K) \to \Gamma$. For $\rho \in S_{K,\Gamma}$,
identify $\rho$ with the linear representation obtained from $\rho$ by
embedding $\Gamma$ into $\GL_n(\CC)$. Put
\begin{equation} \label{eq:totalmass1}
M(K, \Gamma) = \frac{1}{|\Gamma|} \sum_{\rho \in S_{K,\Gamma}} \frac{1}{q^{c(\rho)}},
\end{equation}
assuming that the sum converges. 
\end{defn}
\begin{remark}
It can be shown that the sum in \eqref{eq:totalmass1} is finite if $K$
has mixed characteristics, and is convergent if $K$ has equal 
characteristics, so the definition always makes sense.
\end{remark}

By Burnside's theorem, we can reformulate this definition as follows.
Let $\Sigma_{K,\Gamma}$ be a set of representatives of the 
isomorphism classes (under conjugation within $\Gamma$)
of continuous homomorphisms $\rho: G_K \to \Gamma$.
For $\rho \in \Sigma_{K,\Gamma}$,
let $w(\rho)$ be the order of the centralizer in $\Gamma$
of the image of $\rho$.
Then also 
\begin{equation} \label{eq:totalmass2}
M(K, \Gamma) = \sum_{\rho \in \Sigma_{K,\Gamma}} \frac{1}{w(\rho) q^{c(\rho)}}.
\end{equation}

\begin{remark} \label{rem:summand}
If we view $\Gamma$ as a group equipped with a faithful (i.e., injective)
linear
representation, then adding a trivial summand to the representation
does not change any conductors, and so does not change the total mass.
\end{remark}

\begin{remark}
In some cases, it might be useful to allow $\Gamma$ to be equipped with
a non-faithful linear representation, rather than to view it as a subgroup
of $\GL_n(\CC)$. As this does not materially enrich the situation from our
point of view, we will not do so.
\end{remark}

\begin{lemma} \label{lem:add}
Let $\Gamma_1 \subseteq \GL_m(\CC)$ and $\Gamma_2 \subseteq \GL_n(\CC)$
be finite subgroups, and view $\Gamma_1 \times \Gamma_2$ as a subgroup
of $\GL_{m+n}(\CC)$. Then for any local field $K$,
\[
M(K, \Gamma_1 \times \Gamma_2) = M(K, \Gamma_1) M(K, \Gamma_2).
\]
\end{lemma}
\begin{proof}
Given $\rho: G_K \to \Gamma_1 \times \Gamma_2$ continuous,
let $\rho_1: G_K \to \Gamma_1$ and $\rho_2:
G_K \to \Gamma_2$ be the results of composing $\rho$
with the projections from $\Gamma_1 \times \Gamma_2$ to its two factors.
Then $c(\rho) = c(\rho_1) + c(\rho_2)$, from which the desired result
follows.
\end{proof}

\begin{remark} \label{remark:global}
Bhargava \cite[\S 8]{bhargava} suggests introducing a parameter $s$
in the exponent of $q$ in the definition of the total mass; this 
gives rise to local factors which one then multiplies together
to give a global Dirichlet series, whose asymptotics one hopes
resemble those of a Dirichlet series which actually counts certain
representations of a global Galois group into $\Gamma$. As
verified by Wood \cite{wood}, this heuristic
in fact reproduces Malle's predicted asymptotics for counting
number fields with prescribed Galois group \cite{malle1, malle2}.
We will omit any further consideration in this direction in this paper;
to do so, we omit the parameter $s$, which amounts
to setting $s=1$. (Taking $s$ to be a positive integer amounts
to replacing the linear representation of $\Gamma$ by its $s$-th tensor
power.)
\end{remark}

\section{Bhargava's formula and permutation representations}
\label{sec:perm rep}

Before proving Bhargava's theorem, we first check that the
left-hand side of \eqref{eq:bhargava} is equal to $M(K, S_n)$, where
$S_n$ is embedded in $\GL_n(\CC)$ via its standard permutation representation,
by matching up corresponding terms.
Using this equality, we will establish \eqref{eq:bhargava} in the next
section.

We start with the usual equivalence of categories between
\'etale $K$-algebras and finite $G_K$-sets.
\begin{lemma} \label{lem:biject}
For any field $K$, 
there is a natural bijection between isomorphism classes of
\'etale $K$-algebras of degree $n$ and isomorphism classes,
under conjugation within $S_n$,
of continuous homomorphisms $\rho: G_K \to S_n$.
Under this bijection, finite separable field extensions of $K$
correspond to homomorphisms with transitive image.
\end{lemma}
\begin{proof}
A $K$-algebra $L$ of degree $n$ is \'etale if and only if
there exists an isomorphism of $K^{\sep}$-algebras
$L \otimes_K K^{\sep} \cong (K^{\sep})^n$, that is, if
$L \otimes_K K^{\sep}$ contains $n$ minimal idempotents.
In fact, these idempotents all lie in $L \otimes_K F$, for any 
Galois extension $F/K$ containing a copy of each component of $L$.
Now equip $L \otimes_K K^{\sep}$ with the action of 
$G_K$ which is trivial on the first factor and the
usual action on the second factor. The action on minimal idempotents
yields a continuous homomorphism 
$\rho: G_K \to S_n$.

Conversely, given a continuous homomorphism $\rho: G_K \to
S_n$, we obtain an action of $G_K$ on $(K^{\sep})^n$
via permutations. We can construct a full set of invariants under
some finite separable extension of $K^{\sep}$, which obviously must
coincide with $K^{\sep}$ itself. These invariants form
a $K$-subalgebra of $(K^{\sep})^n$ of $K$-dimension $n$, which
by construction is \'etale. The functors just described yield the desired
bijection.
\end{proof}

We next verify that the bijection of Lemma~\ref{lem:biject}
matches up the two automorphism contributions.
\begin{lemma}
Under the bijection of Lemma~\ref{lem:biject},
let $L$ be an \'etale $K$-algebra corresponding to a continuous
homomorphism $\rho: G_K \to S_n$. Then 
$w(L) = w(\rho)$.
\end{lemma}
\begin{proof}
The group $\Aut(L/K)$ is isomorphic to 
the group of automorphisms of $L \otimes_K K^{\sep}$ which
are semilinear for the $K$-action on $K^{\sep}$. This group in turn coincides
with the $G_K$-equivariant
permutations of the set of minimal idempotents of $L \otimes_K K^{\sep}$,
yielding the claim.
\end{proof}

Finally, we have the equality of the discriminant and conductor
contributions; for this, we need the conductor-discriminant formula.
\begin{lemma} \label{L:conductor}
Let $H$ be an open subgroup of $G = G_K$ with fixed
field $L$,
and let $\rho: H \to \GL_n(\CC)$ be a continuous representation, where
$\GL_n(\CC)$ carries the discrete topology. Then
\[
c(\Ind_H^G \rho) = f(L/K) c(\rho) + n c(L/K).
\]
\end{lemma}
\begin{proof}
See \cite[Proposition~VI.6, Corollary~1]{serre}.
\end{proof}

\begin{lemma} \label{lem:match auto}
Under the bijection of Lemma~\ref{lem:biject},
let $L$ be an \'etale $K$-algebra corresponding to a continuous
homomorphism $\rho: G_K \to S_n$. Then 
$c(L/K) = c(\rho)$.
\end{lemma}
\begin{proof}
Since both functions are additive over direct sums, we may reduce to
the case where $L$ is a field and $\rho$ has transitive image. In that case,
let $F$ be the normal closure of $L/K$, and put $G = \Gal(F/K)$ and
$H = \Gal(F/L)$. Then the linear representation derived from $\rho$
is isomorphic to the representation induced from the trivial
representation on $H$. The desired result thus follows from 
Lemma~\ref{L:conductor} applied to the trivial one-dimensional representation.
\end{proof}

Putting together the three lemmas, we see at once that the left side of
\eqref{eq:bhargava} equals $M(K, S_n)$, as claimed.

\section{Reduction to the mass formula}

We now give a proof of Theorem~\ref{thm:abb}
by reduction to Serre's mass formula. This is essentially Bhargava's
proof in \cite{bhargava} except that we use a standard device from
enumerative combinatorics, the Exponential Formula
(for an exposition of which see \cite[Chapter 5]{stanley}), in lieu of
explicitly combining fields into \'etale algebras as in
\cite[Propositions~1--3]{bhargava}.

\begin{proof}[Proof of Theorem~\ref{thm:abb}]
We first reinterpret Serre's formula as in the previous section.
Let $I$ be the inertia subgroup of $G_K$.
We say a continuous homomorphism $\rho: G_K \to S_n$ with
transitive image is \emph{totally ramified} if 
$\rho^{-1}(S_{n-1})$ and $I$ together generate $G_K$.
Then Serre's formula states that
\[
\frac{1}{n!} \sum_{\rho} \frac{1}{q^{c(\rho)}} =
q^{1-n},
\]
where $\rho$ runs over all totally ramified homomorphisms $G_K \to S_n$ with transitive
image.

Let $T_{K,S_n}$ denote the subset of $S_{K,S_n}$ (in the notation
of Definition~\ref{D:total mass}) consisting
of homomorphisms with transitive image. 
Given $\rho \in T_{K,S_n}$,
let $f = f(\rho)$ 
be the index of the image of $\rho^{-1}(S_{n-1})$ in $G_K/I$;
this index necessarily divides $n$.
Let $K_f$ be the unramified extension of $K$ of degree $f$,
and put $G_f = G_{K_f}$. Then the restriction of $\rho$
to $G_f$ splits as a direct sum of $f$ totally ramified representations
which are isomorphic to the conjugates of some representation
$\rho_f: G_f \to S_{n/f}$ by a generator $\sigma$ of $G_K/G_f$.
By Frobenius reciprocity, $\rho$ is isomorphic to the induced
representation $\Ind^{G_K}_{G_f} \rho_f$, and $c(\rho) = f c(\rho_f)$.

Given a choice of $f$ and $\rho_f$, one can reconstruct such a $\rho$
by choosing a partition of $\{1, \dots, n\}$ into $f$ labeled blocks of
size $n/f$, then choosing a bijection between each group and
$\{1, \dots, n/f\}$. However, each $\rho$ is produced $f/r$ times,
where $r$ is the smallest positive integer such that $\rho_f$
is isomorphic to its conjugate by $\sigma^r$.
In fact the same $\rho$ is produced by all $r$ of the conjugates
of $\rho_f$ under $\sigma$, so we need to divide both by $f/r$ and
by $r$ to account for this.
In addition, there is a further overcount by a factor
of $w(\rho_f)$ (the number of automorphisms of $\rho_f$). Therefore
\[
\frac{1}{n!} \sum_{\rho \in T_{K,S_n}: f(\rho) = f} \frac{1}{q^{c(\rho)}}
= \frac{1}{n!} \frac{n!}{f} 
\sum_{\rho_f} \frac{1}{w(\rho_f) (q^f)^{c(\rho_f)}},
\]
where the sum on $\rho_f$ runs over \emph{isomorphism classes} of
totally ramified representations $\rho_f: G_f \to S_{n/f}$. By
Serre's formula, we have
\begin{equation} \label{eq:intermed}
\frac{1}{n!} \sum_{\rho \in T_{K,S_n}: f(\rho) = f} \frac{1}{q^{c(\rho)}}
= \frac{q^{f-n}}{f}.
\end{equation}

Passing from \eqref{eq:intermed} (after summing over $f$)
to the total mass amounts to
an application of the Exponential Formula:
\[
\begin{split}
\sum_{n=0}^\infty 
M(K, S_n) x^n 
&= \exp\left( \sum_{n=1}^\infty \frac{x^n}{n!} \sum_{\rho \in T_{K,S_n}}
\frac{1}{q^{c(\rho)}} \right) \\
&= \exp \left( \sum_{n=1}^\infty x^n \sum_{f | n} \frac{q^{f-n}}{f}
\right) \\
&= \exp \left( \sum_{i=1}^\infty \sum_{f=1}^\infty \frac{x^{fi} q^{f(1-i)}}{f}
\right) \\
&= \exp \left( \sum_{i=1}^\infty \log (1-x^i q^{1-i})^{-1}
\right) \\
&= \prod_{i=1}^\infty (1-x^i q^{1-i})^{-1}.
\end{split}
\]
Substituting $xq$ for $x$ yields
\[
\sum_{n=0}^\infty 
M(K, S_n) x^n q^n = \prod_{i=1}^\infty (1- x^i q)^{-1}.
\]
The coefficient of $x^n q^{n-k}$ on the right side is visibly equal to
$P(n,n-k)$. Since we checked in Section~\ref{sec:perm rep} that
$M(K,S_n)$ equals the left side of \eqref{eq:bhargava}, we may deduce
the desired result.
\end{proof}

\begin{remark}
Bhargava's proof in \cite{bhargava} shows
that one can formulate more precise versions of the mass formula
that sum over a single splitting type for \'etale algebras of degree $n$,
the case of a single component being \eqref{eq:intermed}. 
In fact, such statements can be read off from the Exponential Formula; in particular,
it will be possible to formulate such results in
other cases (e.g., in Theorem~\ref{thm:bn}), though we will not explicitly do so.
\end{remark}

\section{Uniformity for other groups: tame case}

We now propose a context into which Theorem~\ref{thm:abb}
can potentially be generalized.
\begin{defn}
Let $n$ be a positive integer and let $\Gamma$ be a finite subgroup of 
$\GL_n(\CC)$. If $S$ is a class of local fields,
we say $\Gamma$ is \emph{uniform for $S$} if
there exists a polynomial $P(x) \in \ZZ[x]$ such that for any
local field $K \in S$ with residue field $\FF_q$, we have
$M(K, \Gamma) = P(q^{-1})$.
(There is a natural candidate function $P(q^{-1})$, but it is not always
a polynomial; see Proposition~\ref{prop:tame}.)
If $S$ is the class of all local fields, we say $\Gamma$
is \emph{uniform for local fields}.
\end{defn}

\begin{remark} \label{rem:add}
If $\Gamma_i$ is a finite subgroup of $\GL_{n_i}(\CC)$ which
is uniform for local fields for $i=1,2$, then $\Gamma_1 \times \Gamma_2 \subset
\GL_{n_1 + n_2}(\CC)$ is uniform for local fields, by Lemma~\ref{lem:add}.
\end{remark}

By Theorem~\ref{thm:abb}, the group $S_n \subset \GL_n(\CC)$ is 
uniform for local fields; we may then ask what other groups have this
property. We obtain a candidate formula for the total mass
by calculating what happens in the tamely ramified case, 
i.e., when the residue characteristic
of $K$ is coprime to the order of $\Gamma$. The result is a 
a \emph{quasi-polynomial} in $q^{-1}$, 
i.e., a function which agrees with different polynomials
on different residue classes; the failure of this quasi-polynomial to
be a true polynomial (as in Example~\ref{exa:z3z} below) constitutes a simple
obstruction to $\Gamma$ being uniform for local fields.

\begin{prop} \label{prop:tame}
Let $n$ be a positive integer, and let $\Gamma$ be a finite subgroup of 
$\GL_n(\CC)$.
For $g \in \Gamma$, let $e(g)$ denote the number of eigenvalues of $g$
not equal to $1$, and define the quasi-polynomial $P_\Gamma(q^{-1})$ by
\[
P_\Gamma(q^{-1}) = \frac{1}{|\Gamma|}
\sum_{g,h \in \Gamma: hgh^{-1} = g^q} q^{-e(g)}.
\]
Then  for any local field $K$ whose residue field
$\FF_q$ has characteristic prime to $|\Gamma|$, we have
$M(K, \Gamma) = P_\Gamma(q^{-1})$.
\end{prop}
\begin{proof}
Since $K$ has residue characteristic prime to  $|G|$,
any continuous homomorphism $\rho: G_K \to \Gamma$ factors
through the maximal tame quotient of $K$. That quotient is topologically
generated by $x,y$ subject to the relation $yxy^{-1} = x^q$; thus
the homomorphisms are determined by pairs $(g,h)$ as in the proposition,
and the quantity $e(g)$ is precisely the conductor of the corresponding
homomorphism.
\end{proof}

\begin{cor} \label{cor:tame mass}
With notation as in Proposition~\ref{prop:tame}, we have
\[
M(K, \Gamma) = P_{\Gamma}(q^{-1}) = \sum q^{-e(g)},
\]
where the sum runs over a set of representatives of those conjugacy
classes of $\Gamma$ which are stable under the $q$-th power map.
\end{cor}
\begin{proof}
For fixed $g \in \Gamma$, the set of $h$ such that $hgh^{-1} = g^q$ is
empty if $g$ is not conjugate to $g^q$, and otherwise is a left coset 
of the centralizer of $g$. This yields the claim by Burnside's formula.
\end{proof}
\begin{cor} \label{cor:uniform tame}
Let $n$ be a positive integer and let $\Gamma$ be a finite subgroup of 
$\GL_n(\CC)$. Then $\Gamma$ is uniform for local fields of residue
characteristic prime to $|\Gamma|$ if and only if
the character table of $\Gamma$ has rational integral entries.
\end{cor}
\begin{proof}
By Corollary~\ref{cor:tame mass} plus Dirichlet's theorem on primes in
arithmetic progressions, $\Gamma$ is uniform for local fields
of residue characteristic prime to $|\Gamma|$ if and only if
for each $m$ coprime to $|\Gamma|$, each element of $\Gamma$ is conjugate
to its $m$-th power. It is a standard result of representation
theory for finite groups  \cite[Chapter~13, Theorem~29]{serre-rep}
that this condition
is equivalent to the rationality of the entries of the character
table of $\Gamma$. 
\end{proof}

\begin{remark}
In case $\Gamma$ arises from a permutation representation, $e(g)$
coincides with Malle's index function. This is related to the fact that
one can recover Malle's heuristics by considering local masses;
see Remark~\ref{remark:global}.
\end{remark}

\section{Examples and counterexamples}

In this section, we mention some examples that help clarify
the extent to which uniformity for local fields holds.
We start with an illustration of Proposition~\ref{prop:tame}.
\begin{example} \label{exa:z3z}
If $\Gamma = \ZZ/3\ZZ$
embedded into $\GL_1(\CC)$ as the cube roots of unity, and
$K$ is a local field of residue field $\FF_q$ of characteristic $\neq 3$,
then 
\[
M(K, \Gamma) = \begin{cases}
1 + 2q^{-1} & q \equiv 1 \pmod{3} \\
1 & q \equiv 2 \pmod{3},
\end{cases}
\]
so $\Gamma$ cannot be uniform for local fields.
\end{example}
Example~\ref{exa:z3z} illustrates that the conclusion of
Corollary~\ref{cor:uniform tame} imposes a strong restriction
on groups which can be uniform for local fields.
However, the conclusion of Corollary~\ref{cor:uniform tame} does
not give a sufficient condition for uniformity for all local fields;
that is because it
depends only on the group $\Gamma$ and not on its embedding into
$\GL_n(\CC)$. Here is an example to illustrate what can go wrong
when one changes the embedding.
\begin{example} \label{ex:z2}
Let $\Gamma$ be the group $\ZZ/2\ZZ$, viewed as a subgroup of
$\GL_2(\CC)$ via its regular representation. We know that
$\Gamma$ is uniform for local
fields by Theorem~\ref{thm:abb}; let us check a bit of this explicitly.
For any
local field $K$ whose residue field $\FF_q$ has odd residue
characteristic, $M(K, \Gamma) = 
1 + q^{-1}$ by Proposition~\ref{prop:tame}. 

However, for
$K = \QQ_2$, there are eight continuous homomorphisms
$\rho: G_{\QQ_2}\to \Gamma$, which all
factor through $\Gal(L/\QQ_2)$ for $L = \QQ_2(\zeta_3, i, \sqrt{2})$.
We can compute the conductors of these as 
follows.
Apply local class field theory to 
identify the eight homomorphisms with the homomorphisms 
$\QQ_2^*/(\QQ_2^*)^2 \to \ZZ/2\ZZ$. Writing
$\QQ_2^* = 2^{\ZZ} \times\gotho_{\QQ_2}^* \cong
\ZZ \times \gotho_{\QQ_2}^*$, we can identify
$\QQ_2^*/(\QQ_2^*)^2$ with $\ZZ/2\ZZ \times 
\gotho_{\QQ_2}^*/(\gotho_{\QQ_2}^*)^2$. Now
\begin{align*}
\gotho_{\QQ_2}^* &= 1 + 2 \gotho_{\QQ_2} \\
(\gotho_{\QQ_2}^*)^2 &= 1 + 8 \gotho_{\QQ_2}
\end{align*}
and $\rho: G_{\QQ_2} \to \Gamma$
has conductor $0,2,3$ according to whether
the kernel of the corresponding homomorphism
$\QQ_2^*/(\QQ_2^*)^2 \to \ZZ/2\ZZ$ kills $1+2 \gotho_{\QQ_2}$,
kills $1 + 4 \gotho_{\QQ_2}$ but not $1 + 2 \gotho_{\QQ_2}$,
or does not kill $1 + 4 \gotho_{\QQ_2}$.
There are thus 2 representations of conductor 0, 2 of conductor 2,
and 4 of conductor 3, yielding 
\[
M(\QQ_2, \Gamma) = 1 + 2^{-2} + 2 \cdot 2^{-3} = 1 + 2^{-1}.
\]
In fact, one can make  a similar calculation for any local field of residual
characteristic 2, as in \cite[Exemple~2(b)]{serre-mass} (which in turn follows
\cite[Lemma~4.3]{tunnell}).

Now consider the same group $\Gamma = \ZZ/2\ZZ$, but now embedded
into $\GL_4(\CC)$
via \emph{two} copies of its regular representation. 
(This example will appear again in Proposition~\ref{prop:not g2}.)
Then $P_\Gamma(q^{-1}) = 1 + q^{-2}$, but 
\[
M(\QQ_2, \Gamma) = 1 + 2^{-4} + 2 \cdot 2^{-6} \neq 1 + 2^{-2}.
\]
Hence $\Gamma$ is not uniform for local fields.
In fact $\Gamma$ is not even uniform just for local fields of residual
characteristic 2, as we may see by calculating the total mass over
$\FF_2((t))$. In this case, there are 2
representations of conductor 0, and for each positive integer $i$, there
are $2^i$ representations of conductor $2 \cdot 2i$, yielding 
\[
M(\FF_2((t)), \Gamma) = 
1 + \frac{1}{2} \sum_{i=1}^\infty 2^i \cdot 2^{-4i} = 
1 + \frac{1}{14} \neq 1 + 2^{-4} + 2\cdot 2^{-6}.
\]
\end{example}

\section{Weyl groups}

Serre gave two proofs of his mass formula \eqref{eq:serremass}. One is a
direct computation using $p$-adic integration on a suitable space of
Eisenstein polynomials, but we find more suggestive
the other proof, which is a simple application of the Weyl integration formula
on a rank $n$ division algebra over $K$. Since the units of that
division algebra constitute a
twisted form of the group $\GL_n$ over $K$, Serre's second proof
suggests that the group $S_n$
is arising in Theorem~\ref{thm:abb} as the Weyl group of the Lie
algebra $\mathfrak{gl}_n$, and prompts the following question. (Note that
the semisimple case reduces to the simple case by
Lemma~\ref{lem:add}.)

\begin{question} \label{question:weyl}
Let $\Gamma$ be the Weyl group of a (semi)simple Lie algebra over $\CC$,
embedded in the group of linear transformations of the root space. 
Is $\Gamma$ uniform for local fields of all 
residual characteristics, or if not, for which ones?
\end{question}

Since $W(A_n) = S_n$ equipped with its standard representation,
Theorem~\ref{thm:abb} asserts that the Weyl group $W(A_n)$ is
uniform for local fields. In 
the remainder of the paper, we assemble some additional answers
to Question~\ref{question:weyl}. Namely, let $\Gamma$ be a Weyl group,
and let $S$ be a set of prime numbers. Then 
$\Gamma$ is uniform for local fields of 
residual characteristics in $S$ in each of the following cases:
\begin{itemize}
\item
$\Gamma$ is arbitrary and $S$ consists only of primes not dividing
$|G|$ (Proposition~\ref{prop:tame});
\item
$\Gamma = W(A_n)$ and $S$ is arbitrary (Theorem~\ref{thm:abb});
\item
$\Gamma = W(B_n)$ and $S$ is arbitrary (Theorem~\ref{thm:bn});
\item
$\Gamma = W(D_n)$ and $2 \notin S$ (Theorem~\ref{thm:dn}; see Appendix);
\item
$\Gamma = W(G_2)$ and $2 \notin S$ (Proposition~\ref{prop:g2});
\end{itemize}
but \emph{not} in the following cases:
\begin{itemize}
\item
$\Gamma = W(D_4)$ and $S$ properly contains $\{2\}$
(Proposition~\ref{prop:not d4}; see Appendix);
\item
$\Gamma = W(G_2)$ and $2 \in S$
(Proposition~\ref{prop:not g2}).
\end{itemize}

Keeping in mind the exceptional isomorphisms
$D_2 \cong A_1\times A_1$ and $D_3 \cong A_3$, we see that the remaining
cases of Question~\ref{question:weyl} are:
\begin{itemize}
\item
$\Gamma = W(D_n)$ for $n \geq 5$ and $2 \in S$;
\item
$\Gamma = W(E_6), W(E_7), W(E_8), W(F_4)$.
\end{itemize}
Moreover, we lack an interpretation of Question~\ref{question:weyl}
in terms of $p$-adic integration analogous to Serre's Weyl integration
proof of his formula; the
negative results suggest that any such interpretation
may have to be a bit subtle.

\begin{remark}
One may also pose Question~\ref{question:weyl} for
other finite Coxeter
groups. However, Corollary~\ref{cor:uniform tame} implies
that if 
$\Gamma$ is uniform for local fields, then $\Gamma$ has rational character
table. This is true for all Weyl groups but not typically for other
finite Coxeter groups (like dihedral groups).
It may be better to consider only local fields
which are algebras over an appropriate cyclotomic field (over which the
representations of the group are defined); we have not investigated this
possibility in any detail. Something loosely analogous has been observed
in the global context of counting number fields, where a counterexample to
a conjecture of Malle \cite{malle1}, \cite{malle2} has been given by
Kl\"uners \cite{kluners}, by distinguishing based on the presence or
absence of an appropriate cyclotomic subextension.
\end{remark}

\section{The groups $W(B_n)$}

We now treat Question~\ref{question:weyl} for the Weyl groups
$W(B_n) = W(C_n)$, in a fashion parallel to
that of Section~\ref{sec:perm rep}. Recall that $W(B_n)$ can be identified
with the wreath product of $\ZZ/2\ZZ$ by $S_n$, or the set of $n \times n$
signed permutation matrices; in particular, there is a surjection
$W(B_n) \to S_n$.

\begin{defn}
Given a tower of fields $M/L/K$,
let $w(M/L/K)$ denote the number of automorphisms of the tower 
\emph{over $K$}, that is, preserving but not necessarily fixing $L$.
\end{defn}

\begin{lemma} \label{lem:biject2}
For any field $K$, there is a natural bijection between isomorphism classes
of towers $M/L/K$, where $L/K$ is a separable field extension of degree $n$ and
$M$ is an \'etale $L$-algebra of degree $2$, and isomorphism classes
under conjugation within $W(B_n)$ of continuous homomorphisms
$\rho: G_K \to W(B_n)$ which have transitive image in $S_n$.
\end{lemma}
\begin{proof}
Given a tower $M/L/K$, let $U$ and $V$ be the sets of
minimal idempotents of $L \otimes_K K^{\sep}$
and of $M \otimes_K K^{\sep} = M \otimes_L (L \otimes_K K^{\sep})$,
respectively. Then
$U$ and $V$ form $G_K$-sets of cardinality $n$ and $2n$, respectively,
with $U$ transitive. Moreover, each element of $U$ splits as the sum of
two elements of $V$; this defines a partition of $V$.
We thus obtain a continuous action of $G_K$ on $V$ factoring through $W(B_n)$,
and conversely as in the proof of Lemma~\ref{lem:biject}.
\end{proof}

\begin{lemma} \label{lem:match auto2}
Under the bijection of Lemma~\ref{lem:biject2}, 
$w(M/L/K) = w(\rho)$.
\end{lemma}
\begin{proof}
Analogous to Lemma~\ref{lem:match auto}.
\end{proof}

\begin{lemma} \label{lem:match cond2}
Under the bijection of Lemma~\ref{lem:biject2}, 
$c(\rho) = f(L/K) c(M/L) + c(L/K)$.
\end{lemma}
\begin{proof}
The \'etale algebra $M/L$ corresponds to a permutation representation
of $G_L$ of degree 2, and hence to a one-dimensional linear representation
of $G_L$. The induction of that representation to $G_K$ is precisely the
linear representation corresponding to $\rho$. Hence the claim follows from 
the conductor-discriminant formula (Lemma~\ref{L:conductor}).
\end{proof}

\begin{theorem} \label{thm:bn}
The Weyl group $\Gamma = W(B_n)$ is uniform for local fields.
\end{theorem}
\begin{proof}
We first compute
the contribution to total mass of homomorphisms with transitive
image in $S_n$; let $T_{K,\Gamma}$ be the set of such homomorphisms.
Switching to isomorphism classes, we can rewrite that contribution as a
sum
\[
\sum_\rho \frac{1}{w(\rho) q^{c(\rho)}}
\]
over isomorphism classes (up to conjugation within $W(B_n)$)
of continuous homomorphisms $\rho: G_K \to W(B_n)$ with transitive image in
$S_n$. By Lemmas~\ref{lem:biject2}, \ref{lem:match auto2}, and~\ref{lem:match
cond2}, this sum in turn equals the sum
\[
\sum_{M/L/K} \frac{1}{w(M/L/K) q^{f_{L/K} c(M/L) + c(L/K)}}
= \sum_{M/L/K} \frac{1}{w(M/L/K) q_L^{c(M/L)} q^{c(L/K)}},
\]
where the sum runs over isomorphism classes of towers $M/L/K$ as in
Lemma~\ref{lem:biject2}.

The contribution from towers with $M = L \oplus L$ is
\begin{equation} \label{eq:bn-1}
\sum_{L/K} \frac{1}{2 w(L/K) q^{c(L/K)}} = \sum_{f|n} \frac{q^{f-n}}{2f}
\end{equation}
by \eqref{eq:intermed}. For the other towers, $M$ is a field, so we 
may sum separately over $M/L$ and $L/K$. Before doing so, we account
for automorphisms as in the proof of Theorem~\ref{thm:abb}.
Namely, $w(M/L/K)$ is equal to twice the number $r$ of automorphisms
of $L/K$ which extend to $M$; but when we count separately over $M/L$
and $L/K$, the number of times we count the same tower $M/L/K$
is equal to $w(L/K)/r$, so we need to divide by a factor of
$2r w(L/K)/r = 2w(L/K)$.

We now combine the analysis of the previous paragraph,
Serre's formula
in degree 2, and \eqref{eq:intermed}
to obtain a mass contribution of
\begin{align*}
\sum_{f|n} \sum_{L: f(L/K) = f} \frac{1}{w(L/K) q^{c(L/K)}}
\sum_{M/L} \frac{1}{2 q_L^{c(M/L)}} 
&= \sum_{f|n} \sum_{L: f(L/K) = f} \frac{1}{w(L/K) q^{c(L/K)}}
\left( \frac{1}{2} + q^{-f} \right) \\
&= \sum_{f|n} \frac{q^{f-n}}{f} \left( \frac{1}{2} + q^{-f} \right) \\
&= \sum_{f|n} \left( \frac{q^{f-n}}{2f} + \frac{q^{-n}}{f} \right).
\end{align*}
Putting this together with \eqref{eq:bn-1}, we conclude that the mass
contribution from $T_{K,W(B_n)}$ is 
\begin{equation} \label{eq:mass bn}
\sum_{f|n} \frac{q^{f-n} + q^{-n}}{f}.
\end{equation}
In particular, this contribution is a polynomial in $q^{-1}$, as then is the
total mass by the Exponential Formula.
\end{proof}

\begin{remark}
As in Theorem~\ref{thm:abb}, we may use the
Exponential Formula to compute a generating function for the
total mass. We obtain
\begin{equation} \label{eq:bn gen fn}
\sum_{n=0}^\infty M(K, W(B_n)) x^n = 
 \prod_{i=1}^\infty (1-x^i q^{-i})^{-1} (1-x^i q^{1-i})^{-1}.
\end{equation}
As computed in Theorem~\ref{thm:abb},
\[
\prod_{i=1}^\infty (1-x^i q^{1-i})^{-1} = \sum_{n=0}^\infty
x^n \sum_{k=0}^n \frac{P(n,n-k)}{q^k}.
\]
On the other hand,
\[
\prod_{i=1}^\infty (1-x^i q^{-i})^{-1} = \sum_{n=0}^\infty
x^n q^{-n} P(n),
\]
where $P(n)$ denotes the number of partitions of $n$ into any number
of parts. Hence
the total mass is
\[
q^{-n} \sum_{j=0}^n \sum_{k=0}^j P(j,j-k) P(n-j) q^{j-k}
= q^{-n} \sum_{j=0}^n \sum_{i=0}^j P(j,i) P(n-j) q^i.
\]
\end{remark}

\begin{remark}
It is also possible to check uniformity for
$\Gamma = W(D_n)$ in this fashion, by determining which towers in
Lemma~\ref{lem:biject2} correspond to homomorphisms to $W(B_n)$
whose images lie in $W(D_n)$; this is done in the Appendix
(Theorem~\ref{thm:dn}).
\end{remark}

\begin{remark}
Wood \cite{wood} has generalized Theorem~\ref{thm:bn} to
arbitrary iterated wreath products of symmetric groups, where 
$W(B_n)$ is viewed as the wreath product of $S_2$ by $S_n$. However,
at this level of generality, the mass must be computed using a counting
function which is apparently not the Artin conductor of a linear 
representation; instead, it is a more general function of the images of the
higher ramification groups.
\end{remark}

\section{The group $W(G_2)$}

\begin{convention}
Throughout this section, for $H$ a subgroup of $\Gamma$,
let $\mu(K,H)$ denote the contribution to the
total mass of $(K, \Gamma)$ coming from homomorphisms
$\rho: G_K \to \Gamma$ with image equal to $H$.
\end{convention}

\begin{prop} \label{prop:g2}
The Weyl group $\Gamma = W(G_2)$ is uniform for local fields
of odd residual characteristic.
\end{prop}
\begin{proof}
Identify $W(G_2)$ with the dihedral group 
\[
\Di_6 = \langle s, r | s^6 = r^2 = rsr^{-1} = e \rangle
\]
equipped with its natural two-dimensional representation.
(We write $\Di_n$ to avoid confusion with the Lie algebra $D_n$.)
By Proposition~\ref{prop:tame}, it suffices to check that the total mass of
$(K, \Gamma)$ equals $1 + 2 q^{-1} + 3q^{-2}$ whenever $K$ has residual 
characteristic 3.
We enumerate the subgroups of $\Gamma$ as follows:
\begin{align*}
C_d &= \langle s^{6/d} \rangle \qquad (d=1,2,3,6) \\
D_{d,i} &= \langle rs^i, s^{6/d} \rangle \qquad (d=1,2,3,6; \quad i=0,\dots,6/d - 1).
\end{align*}
Note that $\mu(K, D_{d,i})$ is independent of $i$. 
Write  $\Gamma = D_{3,0} \times C_2$;
given a homomorphism $\rho: G_K \to \Gamma$,
let $\sigma: G_K \to D_{3,0}$ and $\tau: G_K \to \{\pm 1\}$ be the 
homomorphisms induced by the projections.
For $H$ a subgroup of $D_{3,0}$, write $\nu(K, H)$ for the contribution to the
total mass of $(K, \Gamma)$ coming from those $\rho$ for which the associated
$\sigma$ has image $H$.

By applying Theorem~\ref{thm:abb}, we compute
\begin{align*}
\mu(K,C_1) &= \frac{1}{12} \\
\mu(K,C_1) + \mu(K,D_{1,0}) &= \frac{1}{6} + \frac{1}{6} q^{-1} \\
\mu(K, C_1) + 2\mu(K, D_{1,0}) + \mu(K,C_2) + \mu(K, D_{2,0}) 
&= \frac{1}{3} + \frac{2}{3} q^{-1} + \frac{1}{3} q^{-2} \\
\mu(K,C_1) + 3\mu(K,D_{1,0}) + \mu(K,C_3) + \mu(K,D_{3,0}) 
&= \frac{1}{2} + \frac{1}{2} q^{-1} + \frac{1}{2} q^{-2}
\end{align*}
and hence
\begin{align} 
\label{eq:g2-res3-0}
\mu(K, D_{1,0}) &= \frac{1}{12} + \frac{1}{6} q^{-1} \\
\label{eq:g2-res3-0.5}
\mu(K, C_2) + \mu(K, D_{2,0}) &= \frac{1}{12} + \frac{1}{3} q^{-1} + 
\frac{1}{3} q^{-2} \\ 
\mu(K, C_3) + \mu(K, D_{3,0}) &= \frac{1}{6} + \frac{1}{2} q^{-2}.
\label{eq:g2-res3-1}
\end{align}

Remember that we have assumed that $K$ has residual characteristic 3.
By direct calculation, we have
\begin{align}
\nu(K, C_1) &= \frac{2}{12}+ \frac{2}{12}q^{-2} \label{eq:g2-res3-2} \\
\nu(K, D_{1,i}) &= \frac{2}{12} + \frac{8}{12} q^{-1} + \frac{2}{12} q^{-2}
\qquad (i=0,1,2).
\label{eq:g2-res3-3}
\end{align}
Since $\tau$ is at most tamely ramified, we have $c(\rho) = c(\sigma)$ whenever
the linear representation $\sigma$ contains no tamely ramified subrepresentations.
If $\sigma$ has image $C_3$, then it is tamely ramified if and only if it is
unramified. Hence
\begin{equation}
\nu(K, C_3) = \frac{2}{12}(2 + 2 q^{-2}) + 4\left(\mu(K, C_3) - \frac{2}{12}
\right) =
 - \frac{4}{12} + \frac{4}{12} q^{-2} + 4 \mu(K, C_3). \label{eq:g2-res3-4}
\end{equation}
If $\sigma$ instead has image $D_{3,0}$, then it cannot be unramified since
$D_{3,0}$ is not cyclic. It also cannot be tamely ramified: otherwise, if
$M$ were the fixed field of $\ker(\sigma)$, then $M$ would have a quadratic
subextension $L$ over which it would be tame of degree 3, hence unramified,
but then $M/K$ would have an unramified, hence Galois, 
subextension of degree 3 over $K$, and so could not
have Galois group $D_{3,0}$. Thus we have
\begin{equation}
\nu(K, D_{3,0}) = 4\mu(K, D_{3,0}). \label{eq:g2-res3-5}
\end{equation}
Combining \eqref{eq:g2-res3-1} through \eqref{eq:g2-res3-5}, we obtain
\[
\mu(K, \Gamma) = 1 + 2q^{-1} + 3q^{-2},
\]
verifying the desired result.
\end{proof}

\begin{prop} \label{prop:not g2}
The Weyl group $\Gamma = W(G_2)$ is \emph{not} uniform for local fields
of residual characteristic $2$.
\end{prop}
\begin{proof}
Retain notation as in Proposition~\ref{prop:g2}.
In case $K$ has residual characteristic 2,
using \eqref{eq:g2-res3-0} and \eqref{eq:g2-res3-0.5} we obtain
\begin{align}
\nu(K, C_1) &= \frac{1}{12} + \mu(K, C_2) \label{eq:g2-res2-1} \\
\nu(K, D_{1,i}) &= \frac{1}{4} + \frac{2}{3} q^{-1} + \frac{1}{3} q^{-2} - \mu(K, C_2) 
\qquad (i=0,1,2). \label{eq:g2-res2-2}
\end{align}
In this setting, if $\sigma$ has image $C_3$, then it is at most
tamely ramified. Moreover, if $\sigma$ has image $D_{3,0}$, then it is also
at most tamely ramified, as otherwise $D_{3,0}$ would have to have a nontrivial
normal subgroup of 2-power order. In both of these cases,
when $\tau$ is unramified, we have $c(\rho) = c(\sigma)$, and otherwise
we have $c(\rho) = c(\tau)$. Note that there are two tamely ramified
homomorphisms into $C_2$, and there are eight nontrivial homomorphisms
into $C_3$ and $D_{3,0}$: if the residue field of $K$ contains $\FF_4$, 
there are eight homomorphisms into $C_3$ and none into $D_{3,0}$,
and otherwise there are two homomorphisms into $C_3$ and six into $D_{3,0}$.
In either case, we have
\begin{equation}
\nu(K, C_3) + \nu(K, D_{3,0}) = \frac{1}{3} +  q^{-2} + 
8 \left( \mu(K, C_2) - \frac{1}{12} \right). \label{eq:g2-res2-3}
\end{equation}
Adding up \eqref{eq:g2-res2-1}, \eqref{eq:g2-res2-2}, \eqref{eq:g2-res2-3}
yields 
\[
M(K, W(G_2)) = \frac{1}{2} + 2q^{-1} + 2q^{-2} + 6 \mu(K, C_2).
\]
We deduce that $W(G_2)$ is uniform for local fields if and only if 
$C_2$ is; however, we have seen a failure of this in Example~\ref{ex:z2}.
\end{proof}

\begin{remark}
It may be helpful to see how the calculation of total mass in the $G_2$
case works over $K = \QQ_2$, by going through the contributions from
different Galois extensions. (All assertions below may be confirmed
either by direct verification or by consulting the Database of Local Fields 
\cite{data}.)
The trivial extension contributes a mass of
\[
\frac{1}{12}.
\] 
The quadratic extensions were enumerated in Example~\ref{ex:z2};
each quadratic extension contributes one homomorphism with image in
$C_2$ and six with images in the various $D_{0,i}$, for a mass contribution
of 
\[
\frac{1}{12}(7 + 12 \cdot 2^{-2} + 24 \cdot 2^{-3} + 2 \cdot 2^{-4} + 4
\cdot 2^{-6}).
\]
Each $\ZZ/3\ZZ$-extension contributes two homomorphisms; the only
such extension is unramified, 
for a total mass of
\[
\frac{2}{12}.
\]
Each $\ZZ/6\ZZ$-extension contributes two homomorphisms; there are seven
of these, given by the composita of the unramified $\ZZ/3\ZZ$-extension
with each of the seven quadratic extensions. The conductors coincide with
the squares of the conductors of the 
quadratic extensions, yielding a mass contribution of
\[
\frac{1}{12}(2 + 4 \cdot 2^{-4} + 8 \cdot 2^{-6}).
\]
Any homomorphism with image $(\ZZ/2\ZZ)^2$ can be formed from two
homomorphisms with image $\ZZ/2\ZZ$, and the conductor is the sum of
the conductors of those extensions; each unordered 
pair gives six homomorphisms. This yields a mass contribution of
\[
\frac{1}{12}(12 \cdot 2^{-2} + 24 \cdot 2^{-3} + 6 \cdot 2^{-4}
+ 48 \cdot 2^{-5} + 36 \cdot 2^{-6}).
\]
Each $S_3$-extension contributes six homomorphisms; the only such extension
is the Galois closure of
the cubic extension $\QQ_2[z]/(z^3-2)$ of discriminant exponent 2,
yielding a mass contribution of 
\[
\frac{1}{12}(6 \cdot 2^{-2}).
\]
Each $\Di_6 = (S_3 \times C_2)$-extension 
contributes six homomorphisms; there are seven
of these, given by the composita of the $S_3$-extension with the seven
quadratic extensions. The unramified quadratic contributes
$\frac{1}{12}(6 \cdot 2^{-2})$; the other quadratics dominate the conductor,
yielding a contribution from the $\Di_6$-extensions of
\[
\frac{1}{12}(6 \cdot 2^{-2} + 12 \cdot 2^{-4} + 24 \cdot 2^{-6}).
\]
Adding it up yields
\[
M(K, W(G_2)) =
\frac{1}{12}(12 + 36 \cdot 2^{-2} + 48 \cdot 2^{-3} + 24 \cdot 2^{-4} 
+ 48 \cdot 2^{-5} + 72 \cdot 2^{-6}) = \frac{83}{32}.
\]
The calculation of Proposition~\ref{prop:not g2} predicts a total mass of
\[
\frac{1}{2} + 2 \cdot 2^{-1} + 2 \cdot 2^{-2} + 6 \mu(K, C_2),
\]
and from Example~\ref{ex:z2} (after taking out the contribution from
the trivial homomorphism), we have
\[
6 \mu(K, C_2) = \frac{1}{2} + 2^{-4} + 2^{-6}.
\]
So the predicted total mass is also $83/32$, agreeing with the direct
calculation.
\end{remark}

\section{Final remarks}

One may interpret what we have been doing as counting local Galois
representations which are Hodge-Tate with all weights equal to 0. It would
also be natural to try to enumerate $p$-adic Galois representations 
with other Hodge-Tate weights, e.g., two-dimensional
de Rham representations with Hodge-Tate weights 0 and 1. These may 
lead to heuristics for counting global Galois representations with particular
geometric origins, e.g., those arising from the \'etale cohomology of
elliptic curves.

On a more algebro-geometric note, it might make sense to think about counting
representations of $G_K$ into a group $\Gamma$ as counting
$K$-valued points of $B\Gamma$, the classifying stack of $\Gamma$-torsors.
This gives a natural interpretation of the automorphism contribution
to total mass; it is entirely possible that the conductor contribution
also has a natural, possibly Arakelov-theoretic interpretation.

\appendix

\section{Appendix (by Daniel Gulotta)}

In this appendix, we study Question~\ref{question:weyl} for the Weyl group
$W(D_n)$.

Let $\sigma_n: W(B_n) \rightarrow C_2$ be the map whose kernel is $W(D_n)$.
Then $\rho: \Gal(K^{\sep}/K) \rightarrow W(B_n)$ should be counted in the total
mass of $W(D_n)$ iff $\sigma_n \circ \rho$ is the trivial homomorphism.

\begin{lemma}
  Let $L/K$ be a finite algebraic extension of fields.  Let $A$ and
  $A'$ be bases for $L/K$.  Let $f$ be a linear transformation that takes
  $A$ to $A'$.  Then $\Delta(A')=(\det f)^2 \Delta(A)$.
\end{lemma}

\begin{proof}
  See \cite[Proposition 12.1.2]{ireland}.
\end{proof}

Since the ratio of the discriminants of any two bases of $L/K$ is a square,
we can make the following definition.
\begin{defn}
  Let $L/K$ be a finite algebraic extension of fields.  If $\Delta$ is the
  discriminant of some basis of $L/K$, then the \emph{discriminant root field}
  of $L/K$ is $K(\sqrt{\Delta})$. (The terminology is from
\cite{data}.)
\end{defn}

\begin{lemma} \label{L:disctower}
  Let $M/L/K$ be a tower of finite separable algebraic field extensions.
  Let $A$ be a basis for $L/K$, and let $B$ be a basis for $M/L$.
  Let $C=\{ab | a \in A, b \in B\}$.  If $C$ is considered as a basis for
  $M/K$, then $\Delta(C) = \Delta(A)^{[M:L]} \Norm_{L/K} \Delta(B)$.
\end{lemma}

\begin{proof}
  See \cite[Theorem 39]{hilbert}.
\end{proof}

\begin{lemma}
  Let $M/L/K$ be a tower of fields corresponding to a map
  $\rho: \Gal(K^{\sep}/K) \rightarrow W(B_n)$.  Assume that the characteristic
  of $K$ is not $2$.  Then the fixed field of the
  kernel of $\sigma_n \circ \rho$ is the discriminant root field of $M/K$.
\end{lemma}

\begin{proof}
  Apply Lemma \ref{L:disctower}.  Since $[M:L]=2$, $\Delta(C)$ is in the
  same class of $K^*/(K^*)^2$ as $\Norm_{L/K} \Delta(B)$.  For each conjugate
  of $\Delta(B)$, choose one of its square roots.  The product of all of
  these is a square root of $\Norm_{L/K} \Delta(B)$.  Automorphisms
  in $W(D_n)$ will change an even number of signs and will therefore
  preserve the square roots of $\Norm_{L/K} \Delta(B)$, while automorphisms
  in $W(B_n) \setminus W(D_n)$ will flip them.  Hence an element of
  $\Gal(K^{\sep}/K)$ fixes
  the discriminant root field of $M/K$ iff it is in the kernel of
  $\sigma_n \circ \rho$.
\end{proof}

In the following three results, assume that $q$ is odd,
let $\alpha$ be a primitive $(q-1)$-st root of unity, and let
$\pi$ be a uniformizing element of $K$.

\begin{lemma} \label{unram}
If $M/L$ is unramified and $e(L/K)$ is odd, then the discriminant root field
of $M/K$ is $K(\sqrt{\alpha})$.  If $M/L$ is unramified
and $e(L/K)$ is even, then the discriminant root field of $M/K$ is $K$.
\end{lemma}

\begin{proof}
{}From Lemma~\ref{L:disctower}, one can see that the discriminant root field
lies in an unramified extension of $K$, and therefore must be in
$K(\sqrt{\alpha})$.  The Frobenius automorphism acts on the conjugates of
a primitive $(q^{2f}-1)$-st root of unity as a $2f$-cycle.  The total number
of $2f$-cycles is $e$ and each one contributes an odd number of minus signs,
so the Frobenius homomorphism is in $W(D_n)$
iff $e$ is even.  If $n$ is odd, then the discriminant root field cannot be
$K$, so it must be $K(\sqrt{\alpha})$.  If $n$ is even, then the
discriminant root field cannot be $K(\sqrt{\alpha})$ since the Frobenius
automorphism does not fix this field.  Therefore it must be $K$.
\end{proof}

\begin{lemma} \label{ram}
If $M/L$ is ramified and $f(L/K)$ is odd, then the discriminant root field of
$M/K$ is $K(\sqrt{\pi})$ or $K(\sqrt{\alpha \pi})$.
If $M/L$ is ramified and $f(L/K)$ is even, then the
discriminant root field of $M/K$ is $K$ or $K(\sqrt{\alpha})$.  In either
case, the mass is divided equally between the two possibilities.
\end{lemma}

\begin{proof}
Since $M/L$ is ramified of degree two, the discriminant of $M/L$ has
odd valuation in $L$.  Therefore the norm of this element has odd valuation
in $K$ iff $f$ is even.

Let $\beta$ be a primitive $(q^f-1)$-st root of unity.  Then for any
uniformizing
element $\pi_L$ of $L$, Lemma \ref{L:disctower} shows that
the discriminant of $L(\sqrt{\beta \pi_L})/K$ with respect to
the basis $(1,\sqrt{\beta \pi_L})$ is equal to a primitive $(q-1)$st root
of unity times the discriminant of $L(\sqrt{\pi_L})/K$ with respect
to basis $(1,\sqrt{\pi_L})$.  Thus the two extensions have different
discriminant root fields.  They have equal amounts of mass.
\end{proof}

\begin{theorem} \label{thm:dn}
The group $W(D_n)$ is uniform for local fields of odd residual characteristic.
\end{theorem}

\begin{proof}
  Let $G$ be the group of continuous homomorphisms
  $\Gal(K^{\sep}/K) \rightarrow C_2$; this group is abelian.
  For any $\chi \in \hat{G}$, let
  \begin{eqnarray}
    f(\chi) & = &\sum_n \frac{x^n}{2^n n!} \sum_{\rho}
  q^{-c(\rho)} \chi(\sigma_n \circ \rho) \\ \nonumber
  & = & \exp \left[ \sum_n \frac{x^n}{2^n n!}
  \sum_{\rho \mathrm{\ transitive}}
  q^{-c(\rho)} \chi(\sigma_n \circ \rho) \right],
  \end{eqnarray}
where in both cases $\rho$ runs over continuous homomorphisms
from $\Gal(K^{\sep}/K)$ to $W(B_n)$.
  Then the generating function for $M(K, W(D_n))$ is
\[
\sum_{n=0}^\infty M(K, W(D_n)) x^n = 
\frac{4}{|G|}\sum_{\chi \in \hat{G}} f(\chi).
\]
  For $q$ odd, $G$ is isomorphic to $V_4$.  
Define $a$ by
  \[
\log(a)
= \sum_{n=1}^{\infty} \sum_{M/L/K, [L:K]=n} \frac{x^n}{w(M/L/K) 
  q^{c(M/L/K)}},
\]
where the sum runs
  over towers where $M$ splits or is a field with discriminant root field 
  $K$.
  Let $\log(b)$ sum over towers with discriminant root field
  $K(\sqrt{\alpha})$, and
  let $\log(c)$ sum over towers with discriminant root field
  $K(\sqrt{\pi})$; note that $\log(c)$ also equals the sum over towers
with discriminant root field $K(\sqrt{\alpha \pi})$.  
Then 
\[
\sum_{n=0}^\infty M(K, W(D_n)) x^n = 
\frac{1}{4} \left( abc^2+ \frac{ab}{c^2} + 2 \frac{a}{b} \right).
\]
We now imitate the proof
of Theorem~\ref{thm:bn}, using Lemmas \ref{unram} and
  \ref{ram} to sort terms. We get
\[
\log(a) = \sum_{n=1}^\infty x^n \sum_{f|n} \frac{q^{f-n}}{2f} + \sum_{n=1}^\infty
 x^{2n}
\sum_{f|n} \frac{q^{f-2n}}{2f} + \frac{1}{2} \sum_{n=1}^\infty x^{2n} \sum_{f|n} 
\frac{q^{-2n}}{2f}
\]
by adding the sum over $M$ split,
the sum over $M/L$ unramified and $e$ even,
and half the sum over $M/L$ ramified and $f$ even.
We get
\[
2\log(c) = \sum_{n=1}^\infty x^n \sum_{f|n} \frac{q^{-n}}{f}
- \sum_{n=1}^\infty x^{2n} \sum_{f|n} 
\frac{q^{-2n}}{2f}
\]
by adding the sum over $M/L$ ramified, then subtracting the 
sum over $M/L$ ramified and $f$ even.
Exponentiating, adding in \eqref{eq:bn gen fn}, and solving, we get
  \begin{eqnarray}
    a & = & \prod_{n=1}^{\infty} (1-q^{1-n} x^n)^{-1/2}
    (1-q^{1-2n} x^{2n})^{-1/2} (1-q^{-2n} x^{2n})^{-1/4} \\
    b & = & \prod_{n=1}^{\infty} (1-q^{2-2n} x^{2n-1})^{-1/2}
    (1-q^{-2n} x^{2n})^{-1/4} \\
    c & = & \prod_{n=1}^{\infty} (1-q^{-2n} x^{2n})^{1/4}
    (1-q^{-n} x^{n})^{-1/2}.
  \end{eqnarray}
  Therefore 
  \begin{eqnarray} \label{oddmass}
\sum_{n=0}^\infty M(K, W(D_n)) x^n  
  & = & \frac{1}{4} \left[ \prod_{n=1}^{\infty} (1-q^{1-n} x^n)^{-1} \right]
  \left[\prod_{n=1}^{\infty} (1-q^{-n} x^n)^{-1}  \right. \\
\nonumber & & \left. + \prod_{n=1}^{\infty}
    (1+ q^{-n} x^{n})^{-1} \right]+ \frac{1}{2}
  \prod_{n=1}^{\infty} (1-q^{1-2n} x^{2n})^{-1},
  \end{eqnarray}
proving the desired uniformity.
\end{proof}

\begin{prop} \label{prop:not d4}
The group $W(D_4)$ is \emph{not} uniform for all local fields.
\end{prop}
\begin{proof}
Put  $K=\mathbb{Q}_2$ and  
\begin{equation} \label{eq:log a}
\log(a) = \sum_{n=1}^{\infty} \sum_{M/L/K, [L:K]=n} \frac{x^n}{w(M/L/K) 
q^{-c(M/L/K)}},
\end{equation}
where the sum runs over towers where $M$ splits or the discriminant root
field of $M/K$ is $\mathbb{Q}_2$. Let
$\log(b), \log(c), \log(d)$ sum over towers with
discriminant root field 
$\mathbb{Q}_2(\sqrt{-3})$, $\mathbb{Q}_2(\sqrt{-1})$, 
$\mathbb{Q}_2(\sqrt{2})$, respectively. Then 
\[
\sum_{n=0}^\infty M(\QQ_2, W(D_n)) x^n = 
\frac{1}{8} \left( abc^2d^4+ \frac{abc^2}{d^4} + 2 \frac{ab}{c^2}
+ 4\frac{a}{b} \right).
\]
We claim that 
\begin{eqnarray}
\label{eq:q2 a}
a & = & \left[\prod_{n=1}^{\infty} (1-2^{1-n} x^n)^{-1/2} \right]
\exp \left[ \frac{45}{128} x^2 + \frac{11}{256} x^3 + \frac{691}{4096} x^4 +
O(x^5) \right] \\
b & = & \exp \left[ \frac{1}{2} x + \frac{45}{128} x^2 + \frac{257}{768} x^3 +
O(x^4)\right] \\
c & = & \exp \left[ \frac{1}{8} x + \frac{7}{128} x^2 + \frac{23}{768} x^3 +
O(x^4) \right] \\
d & = & \exp \left[ \frac{1}{16} x + \frac{1}{64} x^2 + \frac{1}{192} x^3 +
  O(x^4)\right]
\end{eqnarray}
This precision suffices to imply
$M(\mathbb{Q}_2, W(D_4)) = \frac{1611}{1024}$.
This in turn proves the
desired result: if $W(D_4)$ were uniform for local fields,
we would have to have $M(\mathbb{Q}_2, W(D_4)) = 
\frac{51}{32} = \frac{1632}{1024}$ 
in order to agree with
the coefficient of $x^4$ in \eqref{oddmass} for $q=2$.

It remains to explain how $a,b,c,d$ were computed.
The only calculation which is nontrivial to
verify by hand is the coefficient $\frac{691}{4096}$ in \eqref{eq:q2 a},
i.e., the contribution to \eqref{eq:log a} from terms with $n=4$ and $M$
nonsplit. We do this by feeding data from \cite{data}
into two different programs \cite{gk}.
One (called \texttt{lf}) 
is a combination of C++ and Perl scripts, reading the data from the HTML
served by \cite{data}; the other (called \texttt{gap-check}) 
uses GAP \cite{gap}
within SAGE \cite{sage}
to read in a raw data file available from \cite{data}, then uses SAGE
to tabulate the results.

Here are some details about the verification that the sum of the terms
of \eqref{eq:log a} with $n=4$ and $M$ nonsplit is $\frac{691}{4096}$.
In \cite{data}, we find a table of the 1823 isomorphism classes
of fields $M$ of degree 8 over $K = \QQ_2$; 
each entry includes (among other information)
the discriminant root field of $M$, the order of $\Aut(M/K)$, the
discriminant exponent of $M/K$, the Galois group of the normal closure of $M/K$
as a permutation group, and a list of the isomorphism classes of degree 4 fields
which occur among the subfields of $M$. Data about these degree 4 fields
can be looked up in \cite{data} in an analogous table.

For each $M$ in the table with discriminant root field $K$,
we loop over the isomorphism classes of degree 4 subfields.
We check the Galois groups of the two fields to identify one exceptional case
(see below). In all other cases, there is a unique isomorphism class of towers
$M/L/K$ with $M$ as chosen above and $L$ in the chosen isomorphism class,
and $\#\Aut(M/K) = \#\Aut(M/L/K)$. We compute $c(M/L/K) = c(M/K) - c(L/K)$ by the
conductor-discriminant formula (Lemma~\ref{L:conductor}), and obtain 
one of the desired terms of \eqref{eq:log a}. This count is performed
by the program \texttt{lf} of \cite{gk}.

To determine which Galois groups yields exceptions to the above argument,
we use GAP as follows (see the program \texttt{gap-check} in \cite{gk}).
For each $G$ in a set of representatives under conjugacy for the
transitive subgroups of $S_8$ (a precomputed list in GAP), 
we loop over representatives $H$ of subgroups of $G$
up to conjugacy. For each $H$ and each orbit of $H$ of length 2,
we pick an element $t$ of the orbit, then compute the number of conjugates
$g^{-1} H g$ of $H$ which contain
$\Stab_G(t)$, first for $g$ running over $G$, then over $\Norm_G(H)$, then over
$\Norm_G(\Stab_G(t))$. Call these numbers $c_1, c_2, c_3$. If
$c_1 = c_3$, then there is a unique isomorphism class of towers $M/L/K$ 
with $M$ the fixed field of $G$ and $L$ the fixed field of some conjugate
of $H$. If $c_1 = c_2$, then for any such tower, $\#\Aut(M/K) = \#\Aut(M/L/K)$.

The only exceptional cases are found in the case where $G = \Di_4$, and $H$
equal to a non-normal subgroup of order 2.
(Again, we write $\Di$ for dihedral groups to avoid confusion with the $D$ series
of Lie algebras.) In this case, $c_1 = c_3 \neq c_2$, and in fact
$\#\Aut(M/K) = 8$ and $\#\Aut(M/L/K) = 4$. We simply replace the automorphism
contribution by 4 to obtain the desired term of \eqref{eq:log a}.
\end{proof}


\begin{thebibliography}{99}

\bibitem{belabas}
K. Belabas, Param\'etrisation de structures 
alg\'ebriques et densit\'e de discriminants (d'apr\`es Bhargava),
Sem. Bourbaki 2003/2004, \textit{Ast\'erisque} \textbf{299} (2005),
267--299.

\bibitem{bhargava2}
M. Bhargava, The density of discriminants of quartic rings and fields,
\textit{Ann. Math.} \textbf{162} (2005), 1031--1063.

\bibitem{bhargava3}
M. Bhargava, The density of discriminants of quintic rings and fields,
to appear in \textit{Ann. Math.}

\bibitem{bhargava}
M. Bhargava, Mass formulae for extensions of local fields, and
conjectures on the density of number field discriminants, preprint.

\bibitem{davenport-heilbronn}
H. Davenport and H. Heilbronn, On the density of discriminants of cubic
fields II, \textit{Proc. Royal Soc. London Ser. A} \textbf{322}
(1971), 405--420.

\bibitem{gap}
The GAP Group, GAP -- Groups, Algorithms, and Programming,
version 4.4 (2006), \texttt{http://www.gap-system.org}.

\bibitem{gk}
D. Gulotta and K.S. Kedlaya, scripts available at
\texttt{http://math.mit.edu/\~{}kedlaya/papers}.

\bibitem{humphreys}
J.E. Humphreys, \textit{Reflection Groups and Coxeter Groups},
Cambridge Studies in Advanced Mathematics 29, 
Cambridge Univ. Press (Cambridge), 1990.

\bibitem{data}
J.W. Jones and D.P. Roberts, Database of Local Fields,
\texttt{http://math.la.asu.edu/\~{}jj/localfields}.

\bibitem{kluners}
J. Kl\"uners, A counterexample to Malle's conjecture on the asymptotics
of discriminants, \textit{C.R. Math. Acad. Sci. Paris} \textbf{340} (2005),
411--414.

\bibitem{malle1}
G. Malle, On the distribution of Galois groups, \textit{J. Number Theory}
\textbf{92} (2002), 315--322.

\bibitem{malle2}
G. Malle, On the distribution of Galois groups II, \textit{Exp. Math.}
\textbf{13} (2004), 129--135.

\bibitem{serre-rep}
J.-P. Serre,
\textit{Linear Representations of Finite Groups},
Graduate Texts in Mathematics 42,
Springer-Verlag (New York-Heidelberg), 1977.

\bibitem{serre-mass}
J.-P. Serre, Une ``formule de masse'' pour les extensions totalement
ramifi\'ees de degr\'e donn\'e d'un corps local, \textit{C.R.
Acad. Sc. Paris S\'erie A} \textbf{286} (1978), 1031--1036.

\bibitem{serre}
J.-P. Serre, \textit{Local Fields},
Graduate Texts in Mathematics 67,
Springer-Verlag (New York-Berlin), 1979.

\bibitem{stanley}
R.P. Stanley, \textit{Enumerative Combinatorics}, volume 2,
Cambridge Studies in Advanced Math. 62, Cambridge Univ. Press
(Cambridge), 1999.

\bibitem{sage}
W. Stein and D. Joyner, SAGE: System for Algebra and Geometry Experimentation, 
\textit{Comm. Computer Alg.} \textbf{39} (2005), 61--64; SAGE
version 1.3.6.3 (2006)
available at \texttt{http://sage.math.washington.edu/sage/}.

\bibitem{tunnell}
J.B. Tunnell, On the local Langlands conjecture for $\GL(2)$,
\textit{Inv. Math.} \textbf{46} (1978), 179--200.

\bibitem{wood}
M. Wood, in preparation.

\bibitem{ireland}
K. Ireland and D. Rosen,
\textit{A Classical Introduction to Modern Number Theory},
Graduate Texts in Mathematics 84, Springer-Verlag (New York), 1990.

\bibitem{hilbert}
D. Hilbert, \textit{The Theory of Algebraic Number Fields},
Springer-Verlag (New York), 1998.

\end{thebibliography}
\end{document}